\documentclass[16pt]{article}
    \usepackage[utf8]{inputenc}
    \usepackage{graphicx} 
    \usepackage[margin=1in]{geometry}
    
    \usepackage{amsmath,amssymb,amsthm,bbm}
    \usepackage{mathtools}
    \usepackage{mathabx}
    \usepackage{url}
    \usepackage{color}
    \usepackage{tikz}
    \usepackage{soul}

\usepackage{hyperref}

\hypersetup{
    colorlinks=true,
    citecolor=red,
    filecolor=red,
    linkcolor=red,
    urlcolor=black
}

\usepackage[maxbibnames=99]{biblatex}
\usepackage{biblatex}
\theoremstyle{definition} \newtheorem{dfn}{Definition}[section]
\theoremstyle{definition} 
\theoremstyle{definition} \newtheorem*{axi*}{Axiom}
\theoremstyle{definition} \newtheorem{prp}[dfn]{Proposition}
\theoremstyle{definition} \newtheorem{clm}[dfn]{Claim}
\theoremstyle{definition} \newtheorem{lem}[dfn]{Lemma}
\theoremstyle{definition} \newtheorem{thm}[dfn]{Theorem}
\theoremstyle{definition} \newtheorem{cor}[dfn]{Corollary}
\theoremstyle{definition} 
\theoremstyle{definition} 
\theoremstyle{remark} 

\theoremstyle{remark} \newtheorem{rmk}{Remark}

\newcommand \N {\mathbb N}
\newcommand \Z {\mathbb Z}

\newcommand \R {\mathbb R}

\newcommand \La {\Lambda}

\if10     
\usepackage[mathlines]{lineno}
\newcommand*\patchAmsMathEnvironmentForLineno[1]{%
	\expandafter\let\csname old#1\expandafter\endcsname\csname #1\endcsname
	\expandafter\let\csname oldend#1\expandafter\endcsname\csname end#1\endcsname
	\renewenvironment{#1}%
	{\linenomath\csname old#1\endcsname}%
	{\csname oldend#1\endcsname\endlinenomath}}%
\newcommand*\patchBothAmsMathEnvironmentsForLineno[1]{%
	\patchAmsMathEnvironmentForLineno{#1}%
	\patchAmsMathEnvironmentForLineno{#1*}}%
\AtBeginDocument{%
	\patchBothAmsMathEnvironmentsForLineno{equation}%
	\patchBothAmsMathEnvironmentsForLineno{align}%
	\patchBothAmsMathEnvironmentsForLineno{flalign}%
	\patchBothAmsMathEnvironmentsForLineno{alignat}%
	\patchBothAmsMathEnvironmentsForLineno{gather}%
	\patchBothAmsMathEnvironmentsForLineno{multline}%
}
\linenumbers
\fi

\makeatletter
\newcommand{\blfootnote}[1]{%
  \begingroup
  \renewcommand{\thefootnote}{}%
  \renewcommand{\@makefntext}[1]{\noindent ##1}%
  \footnotetext{#1}%
  \endgroup
}
\makeatother

\newcommand{\GK}[1]{{\color{red}{#1}}}

\title{Spectral Tile Direction in the Group  \(\Z_{p^2} \times \Z_{q^2} \times \Z_r\)}

\author{Thomas Fallon, Gergely Kiss, Azita Mayeli, and  G\'abor Somlai}

\date{\today}

\begin{document}
	
\maketitle

\blfootnote{ 
G.K. was supported by Projects No. FK 142993 and STARTING 150576 of the National Research, Development and Innovation Fund of Hungary, the Bolyai J\'anos Research Fellowship of the Hungarian Academy of Sciences, and the ÚNKP-22-5 New National Excellence Program of the Ministry for Culture and Innovation.
}

\blfootnote{A.M. was supported in part by the National Science Foundation grant DMS-2453769,   
AMS--Simons Research Enhancement Grant, and PSC--CUNY Grant 67807-00 56.}

\blfootnote{G.S. was supported by ARC Discovery Project DP250104965 and OTKA STARTING Grant 150576. 
In much of the work carried out as part of the project, G.S. was a Fulbright Fellow at the Graduate Center of the City University of New York and received support from the Rosztoczy Foundation.}

 \begin{abstract}
Let $p$, $q$, and $r$ be distinct primes such that $p^2q^2<r$. We prove that every spectral set in the cyclic group $\mathbb{Z}_{p^2q^2r}$ is a tile. Since the reverse direction is already known, this shows that $\mathbb{Z}_{p^2q^2r}$ is a Fuglede group under this condition. The proof is based on divisibility properties of mask polynomials and on the structure of spectral sets in finite cyclic groups.
\end{abstract}

	\blfootnote{{\bf Keywords}: Fuglede's conjecture; spectral sets; tiling; finite cyclic groups; mask polynomials.}%

    \blfootnote{%
{\bf Mathematics Subject Classification (2020):}
43A70 (primary); 43A40, 05B45, 20K01 (secondary).}
 
\section{Introduction}

Fuglede's conjecture asks whether spectral sets coincide with translational tiles. Let $G$ be a locally compact abelian group endowed with Haar measure $\mu$, and let $S\subseteq G$ be measurable with $0<\mu(S)<\infty$. We say that $S$ is \emph{spectral} if there is a set of characters whose restrictions to $S$ form an orthogonal basis of $L^2(S)$. We say that $S$ is a \emph{tile} if there is a set $T\subseteq G$ such that almost every element of $G$ has a unique representation of the form $s+t$, where $s\in S$ and $t\in T$. In this case, $T$ is called a \emph{tiling complement} of $S$. We call $G$ a \emph{Fuglede group} if both directions of Fuglede's conjecture hold in $G$.

Fuglede \cite{F1974} proved the conjecture in $\R^n$ under a lattice assumption: a set that tiles by a lattice is spectral, and a spectral set with a lattice spectrum is a tile. The conjecture is false in general. Tao \cite{T2004} constructed spectral sets in finite abelian groups that do not tile and used one of these examples to obtain a counterexample in $\R^5$. Counterexamples were later obtained in dimensions $4$ and $3$; see \cite{FMM2006,KMkishalmaz}. 
In a recent preprint,  Zhang constructed counterexamples to both implications in dimension two \cite{Zhang2026Fuglede}. Thus, the one-dimensional case is the only remaining Euclidean case.

The Euclidean, integer, and finite-group versions of the conjecture are closely related. This connection was initiated by Lagarias and Wang \cite{LW97} through the universal spectrum and universal tiling properties and was further developed in \cite{DJ12,FR2006,KMkishalmaz,PedersenWang2001,LS01,Matolcsi2005}. Dutkay and Lai \cite{DutkayLai2014} gave a systematic account of the relation between the continuous and discrete settings. 
More recently, Fu and Song proved that every normalized spectrum of a bounded spectral set on the real line is rational. Combined with the reduction theorem of Dutkay and Lai, this completes the reduction of the one-dimensional Fuglede conjecture to its finite cyclic analogues.
Thus, finite cyclic groups provide a natural setting in which the conjecture can be studied through arithmetic and algebraic methods, while retaining direct connections with the one-dimensional Euclidean problem.

The structure of tiles and spectral sets in $\Z$ and finite cyclic groups is closely related to the factorization of mask polynomials. Coven and Meyerowitz \cite{CM1999} introduced the conditions (T1) and (T2), which are formulated in terms of cyclotomic divisors of the mask polynomial $m_A$ of a finite set $A\subseteq\Z$. They proved that these conditions imply that $A$ tiles $\Z$, and that they hold for tiles whose cardinality has at most two distinct prime divisors. \L aba \GK{\cite{L2002}} proved that sets satisfying (T1) and (T2) are spectral. These results establish the tile-to-spectral direction for cyclic groups of order $p^nq^m$.
Several cases of the spectral-to-tile direction are also known. The group $\Z_{p^n}$ is a Fuglede group.
Malikiosis and Kolountzakis proved Fuglede's conjecture for $\Z_{p^nq}$ \cite{MK17}, and Malikiosis later extended the spectral-to-tile direction to groups $\Z_{p^mq^n}$ when one of the exponents is at most $6$ \cite{MK20}. 
The fourth author of the present paper   \cite{So2019} proved the spectral-to-tile direction for $\Z_{p^2qr}$, and Zhang \cite{Zhang2024} later generalized this result to $\Z_{p^nqr}$, where $p,q,r$ are distinct primes.  
For squarefree cyclic groups, Kiss, Malikiosis, Somlai, and Vizer \cite{KMSV2021} proved the conjecture for $\Z_{pqrs}$. On the tile-to-spectral side, results for squarefree cyclic groups and for groups whose orders involve two prime divisors follow from \cite{CM1999,Shi18}.  More recently, \L aba and Londner proved the Coven--Meyerowitz conjecture for all integer tilings of period $(p_1p_2p_3)^2$ \cite{LaLo2023,LaLo2024}. Their splitting method also establishes the conjecture for periods
$M=p_1^{n_1}p_2^{n_2}p_3^{n_3}$ whenever $p_1>p_2^{n_2-1}p_3^{n_3-1}$, as well as for certain four-prime families \cite{LaLo2024plus}. 
Together with \L aba's theorem, these results yield the tile-to-spectral direction for the corresponding cyclic groups.

In this paper, we consider the cyclic group $\Z_{p^2q^2r}$, where $p$, $q$, and $r$ are distinct primes and $r>p^2q^2$. By the large-prime splitting result of \L aba and Londner and \L aba's theorem \cite{L2002,LaLo2024plus}, the tile-to-spectral direction is already known in this setting under the weaker assumption $r>pq$. We prove the reverse direction, spectral-to-tile direction, under the stronger assumption $r>p^2q^2$. Our main result is the following. 
 
\begin{thm}\label{thm:main}
Let $p$, $q$, and $r$ be distinct prime numbers. Then every spectral set in the cyclic group $\Z_{p^2q^2r}$ is a tile, provided that $p^2q^2<r$. Consequently, $\Z_{p^2q^2r}$ is a Fuglede group for all distinct primes $p,q,r$ satisfying $p^2q^2<r$.
\end{thm}

We use the condition $r>p^2q^2$ to reduce the possible cardinalities of a spectral set. More precisely, unless the tiling conclusion already follows from an inductive argument, we show that $\Phi_r$ divides the mask polynomial of the spectral set and that $r$ divides its cardinality. The proof then reduces to the cases
$$
r\mid\mid |S|,\qquad pr\mid\mid |S|,\qquad qr\mid\mid |S|,\qquad pqr\mid\mid |S|,
$$
where $d\mid\mid |S|$ means that $\gcd(|G|,|S|)=d$.

The remaining argument combines cyclotomic divisibility with reduction modulo a prime. We organize the relevant divisibility alternatives into a system involving $\Phi_p$, $\Phi_{p^2}$, $\Phi_{pr}$, and $\Phi_{p^2r}$ after reduction modulo $q$. The mod-$q$ method then gives arithmetic restrictions on $|S|$. When one of these divisibility relations fails, the absence of the corresponding Fourier zeros restricts how $S$ may intersect cosets of the subgroups of $\Z_{p^2q^2r}$, and these restrictions lead to explicit tiling complements.

The most involved part occurs when $pqr\mid\mid |S|$ and the full system of divisibility relations holds. In this case, we decompose the spectrum into fibers over $\Z_{p^2r}$ and construct associated spectral pairs in
$$
\Z_{p^2r}\times\Z_q\cong\Z_{p^2qr}.
$$
We then apply the result of the fourth author of the present paper for $\Z_{p^2qr}$ \cite{So2019}, which was later generalized by Zhang to $\Z_{p^nqr}$ \cite{Zhang2024}. To pass from the reduced spectral pairs back to the original group, we use a finite cyclic consequence of the universal tiling theorem of Pedersen and Wang \cite{PedersenWang2001}: several sets with the same tiling spectrum admit a common tiling complement. This allows us first to prove that the spectrum tiles $\Z_{p^2q^2r}$ and then to conclude that the original spectral set also tiles.

The paper is organized as follows. In Section \ref{s2}, we introduce the notation and the main preliminary tools, including mask polynomials, cyclotomic divisibility, the mod-$p$ method, and the cube rule. In Section \ref{s3}, we prove Theorem \ref{thm:main}. We first treat spectral sets of large cardinality, then use the condition $r>p^2q^2$ to reduce the problem to sets whose cardinality is divisible by $r$, and finally consider the remaining cases through the system of cyclotomic divisibility relations. In the Appendix, we prove Proposition \ref{prop:big}, which gives a general tiling criterion when $|G|/\gcd(|G|,|S|)$ is a prime power.

\section{Notation and preliminaries}\label{s2}
Let $S$ be a subset of the cyclic group $\Z_n$. The {\it mask polynomial} of $S$ on a cyclic group is defined as $m_S(x)=\sum_{s\in S}x^s$, where this polynomial is considered as an element of the quotient ring $\Z[x]/(x^n-1)$.

The spectral property of a set is usually expressed by the zeros of its Fourier transform. Notice that in the case of a finite abelian group $G$, the dual group $\widehat{G}$ is isomorphic to $G$, so we identify these two groups.

We say that a set $S$ in $G$ is spectral if there is a set of characters $\La\subseteq\widehat{G}$ such that the restrictions of the elements of $\La$ to $S$ form an orthogonal basis of $L^2(S)$. This implies that $|S|=|\La|$, and in this case we say that $(S,\La)$ is a spectral pair. Using the identification of $G$ and $\widehat{G}$, we can also observe that $(\La,S)$ is a spectral pair.

The fact that a character is a Fourier root of the characteristic function of a set $S$ can be expressed in terms of the mask polynomial of $S$. Let $\chi$ be a character of order $m$, where $m\mid n$. Then $\widehat{1_S}(\chi)=0$ if and only if $m_S$ vanishes at a primitive $m$-th root of unity. This is equivalent to $\Phi_m\mid m_S$, where $\Phi_m$ denotes the $m$-th cyclotomic polynomial. Notice that $\Phi_{p^k}\mid m_S$ implies $p\mid |S|$ for every prime $p$ and $k>0$, since $\Phi_{p^k}(1)=p$.

Another useful observation is that every character can be factorized in the following way. Let $\phi$ be a homomorphism from an abelian group $G$ to $\mathbb{C}^*$. Then there is a cyclic group $\Z_m$ such that $\phi=\chi\circ\pi$, where $\pi$ is a surjective homomorphism from $G$ to $\Z_m$ and $\chi$ is a faithful character of $\Z_m$.

In the case $G=\Z_n$, where $m\mid n$, there is a unique subgroup of $G$ of order $n/m$, which is isomorphic to $\Z_{n/m}$. Thus, there is a natural quotient map $\pi\colon\Z_n\to\Z_m\cong\Z_n/\Z_{n/m}$, which we call the {\it projection} onto $\Z_m$.

The following two lemmas provide some inductive tools proved in \cite{KMSV2021}.
\begin{lem}\label{lem:red}
Suppose that $(S,\La)$ is a spectral pair. Then the following statements hold.
\begin{enumerate}
\item Without loss of generality, we may assume that $0\in S$ and $0\in\La$.
\item If $S$ is contained in a subgroup of $G$, then $S$ is a tile by induction.
\item If $\La$ is contained in a subgroup of $G$, then $S$ is a tile by induction.
\end{enumerate}
\end{lem}
\begin{lem}\label{lemfullcoset}
Suppose that $S$ is a spectral set in $G$ and that $p\mid |G|$ and $S$ is a union of $\Z_p$ cosets. Then $S$ tiles $G$.
\end{lem}

Another technique relies on the fact that cyclotomic polynomials simplify after reduction modulo a prime. We express this in the following lemma. The statements are formulated in the polynomial ring $\Z[x]$, but the corresponding statements also hold after reduction modulo $p$ and in the appropriate quotient rings.

We now recall and prove some statements related to what we call the mod-$p$ method.
\begin{lem}\label{lem:mod}
\begin{enumerate}
\item If $p$ is a prime and $m\in\N$ satisfies $p\nmid m$, then
\begin{equation}\label{eqmp}
\Phi_{mp}(x)=\Phi_m(x^p)/\Phi_m(x).
\end{equation}
\item If $n=p^km$, where $p$ is a prime, $k\geq 1$, and $p\nmid m$, then
\begin{equation}\label{eqmpk}
\Phi_n(x)=\Phi_{pm}(x^{p^{k-1}}).
\end{equation}
\item If $n=p^km$, where $p$ is a prime, $k\geq 1$, and $p\nmid m$, then
$\Phi_n\mid m_S$ in $\Z[x]$ implies $\Phi_m\mid m_S$ in $\Z_p[x]$.
\end{enumerate}
\end{lem}
\begin{proof}
The first two statements are well-known identities for cyclotomic polynomials.

For the third statement, first assume that $n=pm$ and $p\nmid m$. From equation \eqref{eqmp}, we obtain
$$
\Phi_{pm}(x)=\Phi_m(x^p)/\Phi_m(x)=\Phi_m(x)^p/\Phi_m(x)=\Phi_m(x)^{p-1}
$$
in $\Z_p[x]$. This implies that $\Phi_{pm}(x)\mid m_S(x)$ in $\Z[x]$ gives $\Phi_m(x)\mid m_S(x)$ in $\Z_p[x]$.

If $n=p^km$, where $k\geq 2$, then by equation \eqref{eqmpk},
$$
\Phi_n(x)=\Phi_{pm}(x^{p^{k-1}}).
$$ 
After reduction modulo $p$, we have
$$
\Phi_{pm}(x^{p^{k-1}})=\Phi_{pm}(x)^{p^{k-1}}
=\Phi_m(x)^{p^{k-1}(p-1)}.
$$
Therefore, $\Phi_{p^km}(x)\mid m_S(x)$ in $\Z[x]$ implies $\Phi_m(x)\mid m_S(x)$ in $\Z_p[x]$.
\end{proof}

The following consequence of Lemma \ref{lem:mod} is what we call the mod-$p$ method.
\begin{lem}\label{obs1}
Let $n=p^km$, where $p$ is a prime, $k\geq 1$, and $p\nmid m$, and let $S$ be a subset of $\Z_n$. Suppose that
\begin{equation}\label{eqalt}
\Phi_d\mid m_S\text{ in }\Z_p[x]\qquad\text{for every }d\mid m\text{ with }d>1.
\end{equation}
Then $|S|=am+bp$ for some nonnegative integers $a,b\in\Z$.
\end{lem}
\begin{proof}
Let $\overline{S}$ be the projection of $S$ onto $\Z_m$. Notice that $\overline{S}$ is a multiset. The analogue of \eqref{eqalt} also holds for $\overline{S}$, that is, $\Phi_d\mid m_{\overline{S}}$ in $\Z_p[x]$ for every $d\mid m$ with $d>1$.

Since $p\nmid m$, the polynomial $x^m-1$ has no repeated roots in $\Z_p[x]$. Therefore, the cyclotomic polynomials $\Phi_d$, where $d\mid m$, are pairwise relatively prime in $\Z_p[x]$. Also,
$$
\prod_{\substack{d\mid m\ d>1}}\Phi_d(x)=\frac{x^m-1}{x-1}=1+x+\cdots+x^{m-1}.
$$
It follows from \eqref{eqalt} that all coefficients of $m_{\overline{S}}$ are equal modulo $p$. Thus, there is an integer $a$ with $0\leq a<p$ such that $\overline{S}(x)\equiv a\pmod p$ for every $x\in\Z_m$.

Since the multiplicities are nonnegative, we may write $\overline{S}(x)=a+pb_x$ for some nonnegative integer $b_x$. Summing over $x\in\Z_m$, we obtain
$$
|S|=|\overline{S}|=am+p\sum_{x\in\Z_m}b_x.
$$
Taking $b=\sum_{x\in\Z_m}b_x$ proves the statement.
\end{proof}

The following is an easy observation that appears in several papers on Fuglede's conjecture for finite abelian groups whose orders have finitely many distinct prime divisors, such as \cite{KMSV2020, MK20}.
\begin{lem}\label{lem:2primedonotdivide}
Assume that $p$ and $q$ are prime divisors of the order of the cyclic group $\Z_n$. Suppose that $0\in S$ and that $S$ is not contained in any proper subgroup of $\Z_n$. Then there are $s_1,s_2\in S$ such that $p\nmid s_1-s_2$ and $q\nmid s_1-s_2$.
\end{lem}
\begin{proof}
Since $S$ is not contained in any proper subgroup, there is $s_1\in S$ such that $p\nmid s_1$. If $q\nmid s_1$, then $s_2=0$ satisfies the conditions. Thus, we may assume that $q\mid s_1$.

Using the same argument, there is $s_2\in S$ such that $q\nmid s_2$. If $p\nmid s_2$, then the pair $s_2,0$ satisfies the conditions. Thus, we may also assume that $p\mid s_2$. It follows that $p\nmid s_1-s_2$ and $q\nmid s_1-s_2$.
\end{proof}

Using the same argument, the following result was obtained in \cite[Lemma 3.6]{KMSV2021}.
\begin{cor}
Let $(S,\La)$ be a spectral pair in $\Z_n$. Assume that $p$ and $q$ are distinct prime divisors of $n$, and that $p^k\mid n$, $p^{k+1}\nmid n$, $q^l\mid n$, and $q^{l+1}\nmid n$. Then $\Phi_{p^kq^lm}\mid m_S$ and $\Phi_{p^kq^lm'}\mid m_{\La}$ for some $m,m'\mid n/(p^kq^l)$, or $S$ is a tile.
\end{cor}

One of the main tools in the investigation of Fuglede's conjecture is the so-called cube rule. This rule and the following concepts were introduced in \cite{KMSV2020}.

Let $n$ be a positive integer and let $\operatorname{rad}(n)=\prod_{p\mid n}p$ be its square-free radical. Let $\mathcal{P}_n$ denote the set of prime divisors of $n$. Let $\chi$ be a faithful character of $\Z_n$, and let $S\colon\Z_n\to\N$ be a multiset on $\Z_n$. If $\chi$ is a root of the Fourier transform $\widehat{1_S}$ of the weighted characteristic function $1_S$ of $S$, then $\Phi_n\mid m_S$. By \cite[Proposition 3.2]{KMSV2020}, this implies that $S$ is an integer coefficient sum of $\Z_{p_i}$ cosets, where $p_i\in\mathcal{P}_n$ are the distinct prime divisors of $n$.

As a consequence, we obtain the cube rule; see \cite[Proposition 3.5]{KMSV2020}. Let $C=\prod_{p_i\in\mathcal{P}_n}A_i$, where each $A_i\subseteq\Z_{p_i}$ has cardinality $2$. Then the alternating sum of the weights of $S$ over the vertices of $C$ is zero. The sign of each weight is determined by the parity of the Hamming distance from a fixed point. Such a set $C$ is called a cube with respect to the cyclic group $\Z_n$, although $C$ is a product of subsets of the groups $\Z_{p_i}$.

More generally, suppose that $n$ is square-free and $m\mid n$. Let $\overline{S}$ be the projection of $S$ onto $\prod_{p_i\in\mathcal{P}_m}\Z_{p_i}$. The following conditions are equivalent to $\Phi_m\mid m_S$; see \cite{KMSV2020, Stei}.
\begin{enumerate}
\item The multiset $\overline{S}$ is an integer coefficient sum of $\Z_{p_i}$ cosets, where $p_i\in\mathcal{P}_m$.
\item For every cube $C=\prod_{p_i\in\mathcal{P}_m}C_i\subseteq\prod_{p_i\in\mathcal{P}_m}\Z_{p_i}$, where each $C_i$ has cardinality $2$, and every $c_0\in C$, we have
\begin{equation}\label{eqcuberule}
\sum_{c\in C}(-1)^{d_H(c_0,c)}S(c)=0.
\end{equation}
Here, $S(c)$ denotes the weight of $c$ in the multiset $\overline{S}$, and $d_H$ denotes the Hamming distance.
\end{enumerate}

Finally, the authors G.K. and G.S. express their sincere gratitude to the Graduate Center for its unwavering academic support, which greatly contributed to the successful completion of this research.

\section{Proof of the main theorem}\label{s3}

To prove Fuglede's conjecture for $\Z_{p^2q^2r}$, where $r>p^2q^2$, it is enough to prove the spectral-to-tile direction, since the tile-to-spectral direction follows from \L aba's theorem \cite{L2002} and the large-prime splitting result of \L aba and Londner \cite{LaLo2024plus}.

Let $S$ be a spectral set. We will use a case-by-case analysis, but we will try to limit the number of different cases that we have to consider.

We distinguish the cases according to the largest divisor of $p^2q^2r$ that divides the size of $S$. We use the following notation:
$$
d\mid\mid |S|\Longleftrightarrow \gcd(|G|,|S|)=d.
$$

\subsection{Large sets}

\begin{itemize}
\item ${\bf p^2q^2r\mid\mid |S|}$. Then $S=G$, and we are done.

\item ${\bf p^2qr\mid\mid |S|}$ or ${\bf pq^2r\mid\mid |S|}$. We will only prove the first case, since the roles of $p$ and $q$ are symmetric.

Suppose that there are at least $q+1$ points of $S$ in the same $\Z_{q^2}$ coset. Then there are two points whose difference has order $q$ and two points whose difference has order $q^2$. Since $(S,\La)$ is a spectral pair, this gives $\Phi_q\mid m_{\La}$ and $\Phi_{q^2}\mid m_{\La}$. Therefore, $q^2\mid |\La|=|S|$, which is a contradiction.

It follows that every $\Z_{q^2}$ coset contains at most $q$ elements of $S$. Since there are $p^2r$ such cosets, we have $|S|\leq p^2qr$. On the other hand, $p^2qr\mid |S|$, and hence $|S|=p^2qr$. Therefore, each $\Z_{q^2}$ coset contains exactly $q$ elements of $S$. The same argument applies to $\La$.

Since $q^2\nmid |S|$, we cannot have both $\Phi_q\mid m_{\La}$ and $\Phi_{q^2}\mid m_{\La}$. By duality, $(\La,S)$ is also a spectral pair.

If $\Phi_q\nmid m_{\La}$, then each $\Z_q$ coset contains at most one point of $S$. Since each $\Z_{q^2}$ coset contains exactly $q$ points of $S$, these points form a complete set of representatives of the $\Z_q$ cosets in that $\Z_{q^2}$ coset. Hence, $S$ is a tile.

If $\Phi_{q^2}\nmid m_{\La}$, then the intersection of each $\Z_{q^2}$ coset with $S$ is contained in one $\Z_q$ coset. Since this intersection has exactly $q$ elements, it is a complete $\Z_q$ coset. Hence, $S$ is again a tile.

\item ${\bf p^2r\mid\mid |S|}$, ${\bf p^2q^2\mid\mid |S|}$, or ${\bf q^2r\mid\mid |S|}$.
The same type of argument can be adapted to each case, so we only consider one of them.

Suppose that $p^2r\mid\mid |S|$, and write $|S|=kp^2r$. We show that $k=1$ and that $S$ is a tile.

Let $\overline{S}$ be the projection of $S$ onto $\Z_{p^2r}$. Since $q\nmid |S|$, we have $\Phi_q\nmid m_{\La}$ and $\Phi_{q^2}\nmid m_{\La}$. Therefore, each $\Z_{q^2}$ coset contains at most one point of $S$. In other words, the projection $\overline{S}\subseteq\Z_{p^2r}$ is a set.

It follows that $|S|\leq p^2r$. Since $p^2r\mid |S|$, we obtain $|S|=p^2r$. Thus, $\overline{S}$ covers $\Z_{p^2r}$ exactly once, and consequently $S$ is a tile.
\end{itemize}

\begin{rmk}\label{remark:1}
A similar argument shows that if $S$ is a spectral set such that $|G|=p^k\gcd(|G|,|S|)$ for some prime $p$ and $k\in\N$, then $S$ is a tile. This is proved in the Appendix; see Proposition \ref{prop:big}.
\end{rmk}

\subsection{Applying arithmetic conditions}

The purpose of this subsection is to introduce a technique that is useful when one of the prime divisors of $|G|$ is very large. This reduces the number of cases that we need to consider.

We learned the following argument from Ruxi Shi during the Fourier Bases 2018 meeting in Crete\footnote{http://fourier.math.uoc.gr/fb18/}.

\begin{lem}\label{lem:shi}
Let $u$ be a prime divisor of $n$ such that $u^2\nmid n$, and let $S$ be a spectral set in $\Z_n$. Assume that for every proper divisor $d$ of $n$, the spectral-to-tile direction of Fuglede's conjecture holds for $\Z_d$. Suppose further that $u\nmid |S|$ and that $\Phi_{mu}\mid m_S$ implies $\Phi_m\mid m_S$ whenever $m\in\N$ is coprime to $u$. Then $S$ is a tile.
\end{lem}

\begin{proof}
We first show that $uS$ is a set. If $us=us'$ for two distinct elements $s,s'\in S$, then $s-s'$ is a nonzero element of the kernel of multiplication by $u$. Since $u$ is prime and $u^2\nmid n$, the difference $s-s'$ has order $u$. Spectrality then implies $\Phi_u\mid m_{\La}$, and consequently $u\mid |\La|=|S|$, which is a contradiction.

Thus, multiplication by $u$ is injective on $S$, and $uS$ is contained in the subgroup $u\Z_n$, which is isomorphic to $\Z_{n/u}$.

We next show that $uS$ is spectral. Let $\lambda,\lambda'\in\La$ be distinct. If the character corresponding to $\lambda-\lambda'$ has order $m$ with $u\nmid m$, then multiplication by $u$ does not change its order. If its order is $mu$, where $\gcd(m,u)=1$, then $\Phi_{mu}\mid m_S$, and the assumption gives $\Phi_m\mid m_S$. Therefore, the restrictions of the characters in $\La$ remain orthogonal on $uS$. Since $|uS|=|S|=|\La|$, the set $uS$ is spectral.

By the inductive assumption, $uS$ has a tiling complement $T$ in $u\Z_n\cong\Z_{n/u}$. Let $K$ be the kernel of the multiplication map $x\mapsto ux$, and choose one preimage in $\Z_n$ for each element of $T$. Denote the set of these preimages by $T_0$. Then $T_0+K$ is a tiling complement of $S$ in $\Z_n$.

Indeed, applying multiplication by $u$ to the sum $S+(T_0+K)$ gives the tiling $uS+T=u\Z_n$. The injectivity of multiplication by $u$ on $S$, together with the uniqueness of the tiling $uS+T$, gives the uniqueness of the lifted representation. Hence, $S$ tiles $\Z_n$.
\end{proof}

From now on, we assume that ${\bf p^2q^2<r}$. Suppose that the condition of Lemma \ref{lem:shi} does not hold for the prime $r$. Otherwise, $S$ is a tile, and we are done.

The following statement was originally proved in a more general context by \L aba and Marshall in \cite[Lemma 6.1 and Corollary 6.4]{LaMa2022}. We use it to reduce the number of possible cases, and we provide a short proof.

\begin{lem}\label{lem:labalondner}
Assume that $r$ is a prime divisor of the order of the cyclic group $\Z_n$. Suppose that there is an $m\mid n$ such that $\gcd(m,r)=1$, $\Phi_{mr}\mid m_S$, and $\Phi_m\nmid m_S$. Then $|S|\geq r$. Moreover, every $\Z_{n/r}$ coset contains at least one point of $S$.
\end{lem}

\begin{proof}
For a divisor $k$ of $n$, let $\overline{S}_k$ denote the projection of $S$ onto $\Z_k$. Our assumptions mean that the cube rule holds for $\overline{S}_{mr}$, but does not hold for all cubes of $\overline{S}_m$.

The group $\Z_{mr}$ is the union of $r$ cosets of its subgroup $\Z_m$. Fix a cube in the $\Z_m$ direction. The cube rule for $\overline{S}_{mr}$, applied also in the $r$ direction, shows that the alternating sum of the weights on this cube is the same in every $\Z_m$ coset. Therefore, for a fixed cube, either the cube rule holds in every $\Z_m$ coset, or it fails in every $\Z_m$ coset.

Since $\Phi_m\nmid m_S$, there is a cube for which the cube rule fails after projection onto $\Z_m$. For this cube, the corresponding alternating sum is nonzero in every $\Z_m$ coset of $\Z_{mr}$. In particular, every $\Z_m$ coset contains at least one point of $\overline{S}_{mr}$. Since there are $r$ such cosets, we obtain $|\overline{S}_{mr}|=|S|\geq r$.

The same argument shows that the projection of $S$ onto the quotient of order $r$ is surjective. Equivalently, every $\Z_{n/r}$ coset contains at least one point of $S$.
\end{proof}

\begin{lem}\label{lem:phi_r mid}
Assume that the assumptions of Theorem \ref{thm:main} hold. Thus, $S\subseteq\Z_{p^2q^2r}$ is spectral, and $p,q,r$ are distinct primes with $p^2q^2<r$. Then either $S$ is a tile, or
\[
\Phi_r\mid m_S,\qquad \Phi_r\mid m_{\La}.
\]
In the latter case, $r\mid |S|=|\La|$. 
\end{lem}

\begin{proof}
Assume that $S$ is not already a tile. We show that $|S|\geq r$ in this case. It is immediate, if $r\mid |S|$, since then
$|S|\geq r$. Suppose therefore that $r\nmid |S|$.
If the condition of Lemma \ref{lem:shi} holds for $r$, then $S$
is a tile, a contradiction. Otherwise, there is a divisor
$m\mid p^2q^2$ such that
\[
\Phi_{mr}\mid m_S
\qquad\text{and}\qquad
\Phi_m\nmid m_S.
\]
Hence Lemma \ref{lem:labalondner} gives $|S|\geq r$.

There are $p^2q^2$ cosets of the subgroup $\Z_r$. Since $|S|\geq r>p^2q^2$, two distinct points of $S$ belong to the same $\Z_r$ coset. Their difference has order $r$, and spectrality gives $\Phi_r\mid m_{\La}$. Therefore, $r\mid |\La|=|S|$.

Since $|\La|=|S|\geq r>p^2q^2$, the same pigeonhole argument applied to $\La$ gives two distinct elements of $\La$ in the same $\Z_r$ coset. Since $(\La,S)$ is also a spectral pair, this implies $\Phi_r\mid m_S$.
\end{proof}

From now on, we assume ${\bf \Phi_r\mid m_S}$.

\begin{rmk}
The inequality $|S|\geq r$ was also proved by \L aba and Marshall. The additional conclusion of Lemma \ref{lem:labalondner} is that the projection of $S$ onto $\Z_r$ is surjective.

As in the whole argument, we can interchange the roles of $S$ and $\La$. We obtain the following corollary.
\end{rmk}

\begin{cor}\label{cor01}
Let $p^2q^2<r$, and let $(S,\La)$ be a spectral pair. Then either
$S$ is a tile, or
\[
\Phi_r\mid m_S
\qquad\text{and}\qquad
\Phi_r\mid m_{\La}.
\]
In the latter case,
\[
r\mid |S|=|\La|.
\]
\end{cor}


\subsection{Remaining cases}

It remains to prove the statement in the cases
$$
r\mid\mid |S|,\qquad pr\mid\mid |S|,\qquad qr\mid\mid |S|,\qquad pqr\mid\mid |S|.
$$
We may also assume that $\Phi_r\mid m_S$. For completeness, we first list the cases that are excluded by this assumption.

\begin{itemize}
\item ${\bf k\mid\mid |S|}$, where
$$
k\in\{p^2q,pq^2,p^2,q^2,pq,p,q,1\}.
$$
By Corollary \ref{cor01}, either $S$ is a tile or $r\mid |S|$. The second alternative contradicts the assumption on $k$.

\item ${\bf r\mid\mid |S|}$.

Suppose first that $\Phi_{p^2q^2r}\mid m_S$. Then the cube rule holds for every cube in every $\Z_{pqr}$ coset of $\Z_{p^2q^2r}$.

We show that $S$ is a union of $\Z_r$ cosets. We could refer directly to \cite[Lemma 3.5]{So2019}, but we include the argument that is needed here.

\begin{clm}\label{clm1}
Let $n=p^kq^lr^m$, and let $S\subseteq\Z_n$ be a spectral set. Suppose that $\Phi_n\mid m_S$, $\Phi_p\nmid m_{\La}$, and $\Phi_q\nmid m_{\La}$. Then $S$ is a union of $\Z_r$ cosets, and hence $S$ is a tile.
\end{clm}

\begin{proof}
We distinguish two cases. If $S$ is a union of $\Z_r$ cosets, then Lemma \ref{lemfullcoset} implies that $S$ is a tile. Thus, suppose that $S$ is not a union of $\Z_r$ cosets. Then there are $x\in S$ and $y\notin S$ such that $x$ and $y$ belong to the same $\Z_r$ coset.

Let $u\neq x$ belong to the same $\Z_p$ coset as $x$, and let $v\neq x$ belong to the same $\Z_q$ coset as $x$. Since $\Phi_p\nmid m_{\La}$ and $\Phi_q\nmid m_{\La}$, we have $u,v\notin S$. Indeed, if $u\in S$, then $x-u$ has order $p$, and spectrality implies $\Phi_p\mid m_{\La}$. The same argument applies to $v$.

Consider the cube whose three edges from $x$ are determined by $y$, $u$, and $v$. The weights at $x,y,u,v$ are $1,0,0,0$, respectively. Let $z$ be the vertex of this cube whose Hamming distance from $x$ is $3$. By the cube rule,
$$
1+\sum_{\substack{c\in C\\ d_H(x,c)=2}}S(c)-S(z)=0.
$$
Since all the weights are either $0$ or $1$, this equality implies that $S(z)=1$ and that all the vertices at Hamming distance $2$ from $x$ have weight zero. Therefore, $z\in S$.

The choices of $u$ and $v$ were arbitrary. Since $p$ and $q$ are distinct primes, at least one of them is larger than $2$. Suppose first that $p>2$. Choose two distinct elements $u_1,u_2\neq x$ in the $\Z_p$ coset containing $x$, and keep $y$ and $v$ fixed. The corresponding opposite vertices $z_1,z_2$ belong to $S$ and differ only in their $p$-coordinate. Hence $z_1-z_2$ has order $p$, which implies $\Phi_p\mid m_{\La}$, a contradiction.

If $p=2$, then $q>2$. We fix $y$ and $u$ and choose two distinct elements $v_1,v_2\neq x$ in the $\Z_q$ coset containing $x$. The corresponding opposite vertices belong to $S$, and their difference has order $q$. This implies $\Phi_q\mid m_{\La}$, which is again a contradiction.

Therefore, $S$ is a union of $\Z_r$ cosets. Lemma \ref{lemfullcoset} now implies that $S$ is a tile.
\end{proof}

Now suppose that $\Phi_{p^2q^2r}\mid m_{\La}$. Repeating the same argument for $\La$, and using the fact that $(\La,S)$ is a spectral pair, we obtain that $\La$ is a union of $\Z_r$ cosets. By Lemma \ref{lemfullcoset}, $\La$ is a tile.

Since $r\mid\mid |\La|$ and the cardinality of a tile divides $|G|$, it follows that $|\La|=r$. Thus, $\La$ is one $\Z_r$ coset. After a translation, we may assume that $\La=\Z_r$. The spectrality of $(S,\La)$ then implies that $S$ contains exactly one representative from each coset of the annihilator of $\Z_r$, which is the subgroup $\Z_{p^2q^2}$. Therefore, $S$ is a tile.

From now on, suppose that
$$
\Phi_{p^2q^2r}\nmid m_S
\qquad\text{and}\qquad
\Phi_{p^2q^2r}\nmid m_{\La}.
$$
If $S$ is contained in a proper subgroup, then Lemma \ref{lem:red}.2 gives the result. Thus, we may assume that $S$ is not contained in a proper subgroup.

Since $\Phi_{p^2q^2r}\nmid m_{\La}$, every two elements $a,b\in S$ satisfy $p\mid a-b$, $q\mid a-b$, or $r\mid a-b$. Indeed, if none of $p,q,r$ divides $a-b$, then $a-b$ has order $p^2q^2r$, and spectrality would imply $\Phi_{p^2q^2r}\mid m_{\La}$.

By Lemma \ref{lem:2primedonotdivide}, there are $x,y\in S$ such that $p\nmid x-y$ and $q\nmid x-y$. Therefore, $r\mid x-y$.

Let $x'\in S\cap(x+\Z_{pqr})$. Since $x-x'\in\Z_{pqr}$, we have $p\mid x-x'$ and $q\mid x-x'$. Hence $p\nmid y-x'$ and $q\nmid y-x'$. By the previous observation, $r\mid y-x'$. Since $r\mid y-x$, it follows that $r\mid x-x'$.

We claim that the projection of $S$ onto $\Z_{p^2q^2}$ is a set. Suppose that two elements $x,x'\in S$ have the same projection. Then $x-x'\in\Z_r\subseteq\Z_{pqr}$. The argument above gives $r\mid x-x'$. On the other hand, every element of the subgroup $\Z_r$ is divisible by $p^2q^2$. Since $p,q,r$ are distinct, this implies $x-x'=0$. Thus, the projection is injective, and $|S|\leq p^2q^2$.

This is impossible, because $r\mid |S|$, $S\neq\emptyset$, and $p^2q^2<r$. Therefore, this case cannot occur unless $S$ is contained in a proper subgroup, and in that case Lemma \ref{lem:red}.2 shows that $S$ is a tile.

\item ${\bf pr\mid\mid |S|}$ or ${\bf qr\mid\mid |S|}$.

Since the roles of $p$ and $q$ are symmetric, it is enough to consider the case ${\bf pr\mid\mid |S|}$.

In this case,
$$
\Phi_q\nmid m_S,\qquad \Phi_{q^2}\nmid m_S,
$$
and
$$
\Phi_q\nmid m_{\La},\qquad \Phi_{q^2}\nmid m_{\La}.
$$

Project $S$ onto $\Z_{p^2r}$. This projection, denoted by $\overline{S}$, is a set. Otherwise, two points of $S$ would belong to the same $\Z_{q^2}$ coset. Their difference would have order $q$ or $q^2$, which would imply $\Phi_q\mid m_{\La}$ or $\Phi_{q^2}\mid m_{\La}$. This is impossible.

Since $\overline{S}$ is a set, we have $|S|\leq p^2r$. Since $p^2r\nmid |S|$, we have $|S|<p^2r$.

Because $pr\mid |S|$, at least one $\Z_{p^2}$ coset in $\Z_{p^2r}$ contains at least $p$ elements of $\overline{S}$.

Assume first that there is a $\Z_{p^2}$ coset containing at least $p+1$ elements of $\overline{S}$. Thus, there is an $x\in\Z_{p^2r}$ such that
$$
|\overline{S}\cap(\Z_{p^2}+x)|\geq p+1.
$$
Since $\Phi_r\mid m_S$, all $\Z_{p^2}$ cosets contain the same number of points of $\overline{S}$. Therefore, each of them contains at least $p+1$ points.

The points of $\overline{S}$ have lifts to $S$ whose differences may contain powers of $q$. The pigeonhole principle gives
$$
\Phi_p\mid m_{\La}
\quad\text{or}\quad
\Phi_{pq}\mid m_{\La}
\quad\text{or}\quad
\Phi_{pq^2}\mid m_{\La},
$$
and also
$$
\Phi_{p^2}\mid m_{\La}
\quad\text{or}\quad
\Phi_{p^2q}\mid m_{\La}
\quad\text{or}\quad
\Phi_{p^2q^2}\mid m_{\La}.
$$
Using points in different $\Z_{p^2}$ cosets, together with $r>p$ and $\Phi_r\mid m_{\overline{S}}$, we similarly obtain
$$
\Phi_{pr}\mid m_{\La}
\quad\text{or}\quad
\Phi_{pqr}\mid m_{\La}
\quad\text{or}\quad
\Phi_{pq^2r}\mid m_{\La},
$$
and
$$
\Phi_{p^2r}\mid m_{\La}
\quad\text{or}\quad
\Phi_{p^2qr}\mid m_{\La}
\quad\text{or}\quad
\Phi_{p^2q^2r}\mid m_{\La}.
$$

By Lemma \ref{lem:mod}, these relations imply
\begin{equation}\label{eqrel1}
\Phi_p\mid m_{\La},\quad
\Phi_{p^2}\mid m_{\La},\quad
\Phi_r\mid m_{\La},\quad
\Phi_{pr}\mid m_{\La},\quad
\Phi_{p^2r}\mid m_{\La}
\end{equation}
in $\Z_q[x]/(x^n-1)$.

If all the divisibility relations in \eqref{eqrel1} hold, then Lemma \ref{obs1}, applied with $m=p^2r$, gives
$$
|S|=ap^2r+bq
$$
for some nonnegative integers $a$ and $b$. Since $|S|<p^2r$, we have $a=0$. Therefore, $q\mid |S|$, which is a contradiction.

It remains to consider the case in which every $\Z_{p^2}$ coset in $\Z_{p^2r}$ contains exactly $p$ elements of $\overline{S}$. Since $\Phi_r\mid m_S$, all these cosets have the same cardinality, and hence $|S|=pr$. We may also assume that at least one of the relations in \eqref{eqrel1} fails in $\Z_q[x]/(x^n-1)$.

If $\Phi_p\nmid m_{\La}$ in $\Z_q[x]$, then none of
$$
\Phi_p,\qquad \Phi_{pq},\qquad \Phi_{pq^2}
$$
divides $m_{\La}$ in $\Z[x]$. Thus, every $\Z_p$ coset of $\Z_{p^2r}$ contains at most one element of $\overline{S}$. Since each $\Z_{p^2}$ coset contains exactly $p$ elements, these elements form a complete set of representatives of the $\Z_p$ cosets. Hence, $S$ is a tile.

If $\Phi_{p^2}\nmid m_{\La}$ in $\Z_q[x]$, then none of
\[
\Phi_{p^2},\qquad \Phi_{p^2q},\qquad \Phi_{p^2q^2}
\]
divides $m_{\La}$ in $\Z[x]$. Therefore, the $p$ points of
$\overline{S}$ in each $\Z_{p^2}$ coset are contained in a
single $\Z_p$ coset. Since each such intersection has exactly
$p$ elements, it is a complete $\Z_p$ coset.

Let $T_p$ be a complete system of representatives of the
$\Z_p$ cosets in $\Z_{p^2}$. Then
\[
\overline{S}+T_p=\Z_{p^2r}
\]
is a tiling. Since the projection of $S$ onto $\overline{S}$
is injective, lifting this tiling and adjoining the kernel
$\Z_{q^2}$ gives a tiling of $G$. Hence, $S$ is a tile.


If $\Phi_{pr}\nmid m_{\La}$ in $\Z_q[x]$, then none of
\[
\Phi_{pr},\qquad \Phi_{pqr},\qquad \Phi_{pq^2r}
\]
divides $m_{\La}$ in $\Z[x]$. It follows that, in each
$\Z_{pr}$ coset, the intersection with $\overline{S}$ is
contained either in a $\Z_p$ coset or in a $\Z_r$ coset.

There are $p$ cosets of $\Z_{pr}$ in $\Z_{p^2r}$. An
intersection of the first type has at most $p$ elements,
whereas an intersection of the second type has at most $r$
elements. Since $|\overline{S}|=pr$ and $r>p$, no intersection
can be of the first type: otherwise,
\[
|\overline{S}|
\leq p+(p-1)r<pr.
\]
Consequently, every $\Z_{pr}$ coset contains exactly one
complete $\Z_r$ coset of $\overline{S}$.

It follows that
\[
\overline{S}+\Z_p=\Z_{p^2r}
\]
is a tiling. Since the projection of $S$ onto $\overline{S}$
is injective, lifting this tiling and adjoining the kernel
$\Z_{q^2}$ gives a tiling of $G$. Hence, $S$ is a tile.


Finally, suppose that $\Phi_{p^2r}\nmid m_{\La}$ in $\Z_q[x]$. If $\La$ is contained in a proper subgroup, then Lemma \ref{lem:red}.3 gives the result. Otherwise, Lemma \ref{lem:2primedonotdivide}, together with the corollary following it, gives
$$
\Phi_{p^2r}\mid m_{\La},
\quad\text{or}\quad
\Phi_{p^2qr}\mid m_{\La},
\quad\text{or}\quad
\Phi_{p^2q^2r}\mid m_{\La}.
$$
In each case, Lemma \ref{lem:mod} gives $\Phi_{p^2r}\mid m_{\La}$ in $\Z_q[x]$, which is a contradiction.

The following lemma will be used in the last case. It is a finite cyclic consequence of the universal tiling result of Pedersen and Wang. Its role is to show that, once the reduced sets have the same spectrum and that spectrum is a tile, their tiling complements can be chosen uniformly. 

\begin{lem}\label{lem:common-tiling-spectrum}
Let $N$ be a positive integer, and let $A_1,\ldots,A_t,\Gamma\subseteq\Z_N$. Suppose that $\Gamma$ is a tile of $\Z_N$ and that $(A_j,\Gamma)$ is a spectral pair for every $j$. Then the sets $A_1,\ldots,A_t$ have a common tiling complement in $\Z_N$.
\end{lem}

\begin{proof}
We identify all the sets with subsets of $\{0,\ldots,N-1\}$. For every $j$, let $\Omega_j=A_j+[0,1)$. Since $(A_j,\Gamma)$ is a spectral pair in $\Z_N$, the set
$$
\frac{1}{N}\Gamma+\Z
$$
is a spectrum of $\Omega_j$.

Since $\Gamma$ tiles $\Z_N$, it is a direct summand of $\Z_N$. Therefore, $\frac{1}{N}\Gamma+\Z$ is a strongly periodic set in the terminology of \cite{PedersenWang2001}. By \cite[Theorem 4.5]{PedersenWang2001}, there is a strongly periodic tiling set, depending only on $\frac{1}{N}\Gamma+\Z$, which tiles every $\Omega_j$. In the construction of that theorem, the tiling set can be written as $T+N\Z$ for some $T\subseteq\{0,\ldots,N-1\}$. Therefore,
$$
\Omega_j+(T+N\Z)=\R
$$
is a tiling for every $j$. Since $\Omega_j$ is a union of unit intervals with integer endpoints, reducing this tiling modulo $N$ gives
$$
A_j+T=\Z_N
$$
with unique representations for every $j$. Hence, $T$ is a common tiling complement of the sets $A_1,\ldots,A_t$.
\end{proof}

\item ${\bf pqr\mid\mid |S|}$.

Suppose that $|S|=upqr$ for some nonzero $u\in\N$ with $\gcd(u,pqr)=1$. We have $\Phi_r\mid m_S$ and $\Phi_r\mid m_{\La}$. The following table describes the main strategy of the proof.

\begin{table}[h!]
\centering
\begin{tabular}{ccc}
$\Phi_p\mid m_{\La}$ or $\Phi_{pq}\mid m_{\La}$ or $\Phi_{pq^2}\mid m_{\La}$
& $\Longrightarrow$ &
$\Phi_p\mid m_{\La}$ in $\Z_q[x]/(x^n-1)$\\
$\Phi_{pr}\mid m_{\La}$ or $\Phi_{pqr}\mid m_{\La}$ or $\Phi_{pq^2r}\mid m_{\La}$
& $\Longrightarrow$ &
$\Phi_{pr}\mid m_{\La}$ in $\Z_q[x]/(x^n-1)$\\
$\Phi_{p^2}\mid m_{\La}$ or $\Phi_{p^2q}\mid m_{\La}$ or $\Phi_{p^2q^2}\mid m_{\La}$
& $\Longrightarrow$ &
$\Phi_{p^2}\mid m_{\La}$ in $\Z_q[x]/(x^n-1)$\\
$\Phi_{p^2r}\mid m_{\La}$ or $\Phi_{p^2qr}\mid m_{\La}$ or $\Phi_{p^2q^2r}\mid m_{\La}$
& $\Longrightarrow$ &
$\Phi_{p^2r}\mid m_{\La}$ in $\Z_q[x]/(x^n-1)$
\end{tabular}
\caption{System of divisibility relations}
\label{tab1}
\end{table}

{\bf Case 1.}
Suppose that in each line of Table \ref{tab1}, at least one of the divisibility relations on the left-hand side holds. By Lemma \ref{lem:mod}, we have
$$
\Phi_r\mid m_{\La},\quad
\Phi_p\mid m_{\La},\quad
\Phi_{pr}\mid m_{\La},\quad
\Phi_{p^2}\mid m_{\La},\quad
\Phi_{p^2r}\mid m_{\La}
$$
in $\Z_q[x]/(x^n-1)$.

Lemma \ref{obs1} gives
$$
|S|=upqr=kp^2r+lq
$$
for some nonnegative integers $k$ and $l$. Therefore,
$$
kp^2r=q(upr-l),
$$
and since $q\nmid p^2r$, we have $q\mid k$.

If $k>0$, then $k\geq q$, and hence $|S|\geq p^2qr$. No $\Z_{q^2}$ coset can contain more than $q$ elements of $S$. Otherwise, two differences of orders $q$ and $q^2$ would occur, giving
$$
\Phi_q\mid m_{\La}
\qquad\text{and}\qquad
\Phi_{q^2}\mid m_{\La}.
$$
This would imply $q^2\mid |\La|=|S|$, which is a contradiction.

Since there are $p^2r$ cosets of $\Z_{q^2}$ and each contains at most $q$ elements of $S$, we have $|S|\leq p^2qr$. Therefore, $|S|=p^2qr$, and every $\Z_{q^2}$ coset contains exactly $q$ elements of $S$.

Since $q^2\nmid |S|$, at least one of $\Phi_q\mid m_{\La}$ and $\Phi_{q^2}\mid m_{\La}$ fails. If $\Phi_q\nmid m_{\La}$, then the $q$ points in each $\Z_{q^2}$ coset form a system of representatives of the $\Z_q$ cosets. If $\Phi_{q^2}\nmid m_{\La}$, then the $q$ points form one complete $\Z_q$ coset. In both cases, $S$ is a tile.

Now suppose that $k=0$. Then $|S|=lq$. In this case, Lemma \ref{obs1} has been applied to $m_{\La}$. The proof of that lemma shows that every multiplicity in the projection of $\La$ onto $\Z_{p^2r}$ is divisible by $q$.

Every $\Z_{q^2}$ coset contains at most $q$ elements of $\La$. Indeed, if one of these cosets contained at least $q+1$ elements, then there would be two differences of orders $q$ and $q^2$. Since $(\La,S)$ is a spectral pair, this would imply $\Phi_q\mid m_S$ and $\Phi_{q^2}\mid m_S$, and hence $q^2\mid |S|$, which is a contradiction. Therefore, every $\Z_{q^2}$ coset contains either zero or exactly $q$ elements of $\La$.

Since $q^2\nmid |S|$, at least one of the relations $\Phi_q\mid m_S$ and $\Phi_{q^2}\mid m_S$ fails. If $\Phi_{q^2}\nmid m_S$, then the $q$ elements of $\La$ in every nonempty $\Z_{q^2}$ coset are contained in one $\Z_q$ coset. Thus, $\La$ is a union of $\Z_q$ cosets, and Lemma \ref{lemfullcoset} shows that $\La$ is a tile. Applying Lemma \ref{lem:common-tiling-spectrum} to the spectral pair $(S,\La)$, we conclude that $S$ is a tile.

It remains to consider the case
$$
\Phi_q\nmid m_S.
$$
Then the $q$ elements of $\La$ in every nonempty $\Z_{q^2}$ coset form a complete system of representatives of the $\Z_q$ cosets.  We now separate the two $q$-adic digits of the $\Z_{q^2}$-coordinate. Thus, in every nonempty fiber of $\La$, the first $q$-adic digit runs through all elements of $\Z_q$, while the second digit may depend on the first one.

Set $M=p^2r$. Since $\gcd(M,q)=1$, we identify
$$
\Z_{p^2q^2r}\cong\Z_M\oplus\Z_{q^2},
$$
and we use the corresponding identification of the dual group. There is a set $A\subseteq\Z_M$ and, for every $a\in A$, a function $f_a\colon\Z_q\to\Z_q$ such that
$$
\La=\{(a,j+qf_a(j)):a\in A,\ j\in\Z_q\}.
$$
Since every nonempty fiber has cardinality $q$, we have
$$
|A|=\frac{|\La|}{q}=upr.
$$

Two distinct elements of $\La$ in the same fiber have a difference of order $q^2$. Since $(\La,S)$ is a spectral pair, it follows that $\Phi_{q^2}\mid m_S$.

For $t,h\in\Z_q$, define
$$
S_{t,h}=\{x\in\Z_M:(x,t+qh)\in S\}.
$$
Let $\omega_d=e^{2\pi i/d}$. The relation $\Phi_{q^2}\mid m_S$ gives
$$
0=\sum_{t=0}^{q-1}\omega_{q^2}^t
\sum_{h=0}^{q-1}|S_{t,h}|\omega_q^h.
$$
The elements $1,\omega_{q^2},\ldots,\omega_{q^2}^{q-1}$ are linearly independent over $\mathbb{Q}(\omega_q)$. Therefore,
$$
\sum_{h=0}^{q-1}|S_{t,h}|\omega_q^h=0
$$
for every $t\in\Z_q$. Since the minimal polynomial of $\omega_q$ is $1+x+\cdots+x^{q-1}$, the number $|S_{t,h}|$ does not depend on $h$. We denote this common number by $n_t$.

For $h,j\in\Z_q$, define
$$
B_h=\{(x,t)\in\Z_M\oplus\Z_q:x\in S_{t,h}\}
$$
and
$$
\Gamma_j=\{(a,f_a(j)):a\in A\}.
$$
We claim that $(B_h,\Gamma_j)$ is a spectral pair in $\Z_M\oplus\Z_q$ for every $h,j\in\Z_q$.

Let $(x,t),(x',t')$ be two distinct elements of $B_h$. Then $(x,t+qh),(x',t'+qh)\in S$.  By duality, $(\La,S)$ is also a spectral pair. Therefore, the orthogonality of the characters indexed by $S$ on $\La$ gives 
$$
0=\sum_{a\in A}\sum_{\ell=0}^{q-1}
\omega_M^{(x-x')a}
\omega_{q^2}^{(t-t')(\ell+qf_a(\ell))}.
$$
If $t=t'$, then
$$
\sum_{a\in A}\omega_M^{(x-x')a}=0,
$$
which is the required orthogonality relation for $\Gamma_j$.

Suppose that $t\neq t'$, and put $d=t-t'$. For $\ell\in\Z_q$, let
$$
C_\ell=\sum_{a\in A}
\omega_M^{(x-x')a}\omega_q^{df_a(\ell)}.
$$
Then
$$
\sum_{\ell=0}^{q-1}C_\ell\omega_{q^2}^{d\ell}=0.
$$
Each $C_\ell$ belongs to $\mathbb{Q}(\omega_{Mq})$. Since $q\nmid M$, we have
$$
[\mathbb{Q}(\omega_{Mq^2}):\mathbb{Q}(\omega_{Mq})]=q.
$$
Also, $\omega_{q^2}^d$ is a primitive $q^2$-th root of unity. Hence,
$$
1,\omega_{q^2}^d,\ldots,\omega_{q^2}^{d(q-1)}
$$
are linearly independent over $\mathbb{Q}(\omega_{Mq})$. It follows that $C_\ell=0$ for every $\ell$. In particular,
$$
\sum_{a\in A}
\omega_M^{(x-x')a}\omega_q^{(t-t')f_a(j)}=0.
$$
Thus, the characters indexed by $\Gamma_j$ are orthogonal on $B_h$.

Moreover,
$$
|B_h|=\sum_{t\in\Z_q}n_t=\frac{|S|}{q}=upr=|A|=|\Gamma_j|.
$$
Therefore, $(B_h,\Gamma_j)$ is a spectral pair.

The group $\Z_M\oplus\Z_q$ is isomorphic to $\Z_{p^2qr}$. By the result of the fourth author  \cite{So2019}, every $B_h$ is a tile of $\mathbb Z_{p^2qr}$.
This result is also a special case of Zhang's later theorem for the
more general family $\mathbb Z_{p^nqr}$ \cite{Zhang2024}. Hence $|B_h|$ divides $p^2qr$. Since $|B_h|=upr$ and $\gcd(u,pqr)=1$, we obtain $u=1$. In particular,
$$
|S|=pqr
\qquad\text{and}\qquad
|A|=pr.
$$

Fix $h=0$. For every $j\in\Z_q$, $(\Gamma_j,B_0)$ is a spectral pair by duality. Since $B_0$ is a tile of $\Z_M\oplus\Z_q$, Lemma \ref{lem:common-tiling-spectrum} gives a set $T\subseteq\Z_M\oplus\Z_q$ such that
$$
\Gamma_j+T=\Z_M\oplus\Z_q
$$
is a tiling for every $j\in\Z_q$.

Let
$$
K=\Z_M\oplus q\Z_{q^2},
$$
which is a subgroup of $\Z_M\oplus\Z_{q^2}$ isomorphic to $\Z_M\oplus\Z_q$. For $j\in\Z_q$, define
$$
\widetilde{\Gamma}_j=\{(a,qf_a(j)):a\in A\}\subseteq K,
$$
and let $\widetilde{T}\subseteq K$ be the image of $T$ under the isomorphism
$$
\Z_M\oplus\Z_q\to K,\qquad (x,h)\mapsto(x,qh).
$$
Then
$$
\widetilde{\Gamma}_j+\widetilde{T}=K
$$
is a tiling for every $j\in\Z_q$. Also,
$$
\La=\bigcup_{j=0}^{q-1}\bigl((0,j)+\widetilde{\Gamma}_j\bigr),
$$
and the sets in this union lie in the different cosets of $K$. Therefore,
$$
\La+\widetilde{T}=\Z_M\oplus\Z_{q^2}
$$
is a tiling. Thus, $\La$ is a tile of $\Z_{p^2q^2r}$.

Finally, applying Lemma \ref{lem:common-tiling-spectrum} to the spectral pair $(S,\La)$ shows that $S$ is a tile.

{\bf Case 2.}
Suppose that there is a line of Table \ref{tab1} for which none of the divisibility relations on the left-hand side holds.

We will frequently use the following observation. Since $p^2\nmid |S|$, we cannot have both $\Phi_p\mid m_{\La}$ and $\Phi_{p^2}\mid m_{\La}$. Therefore, $S$ has at most $p$ elements in every $\Z_{p^2}$ coset. Similarly, since $q^2\nmid |S|$, the set $S$ has at most $q$ elements in every $\Z_{q^2}$ coset.

We distinguish four cases.

\begin{enumerate}
\item\label{item:subcase1}
Suppose that
$$
\Phi_p\nmid m_{\La},\qquad
\Phi_{pq}\nmid m_{\La},\qquad
\Phi_{pq^2}\nmid m_{\La}.
$$
These conditions imply that in each $\Z_{pq^2}$ coset, at most one $\Z_{q^2}$ coset can contain elements of $S$. Since every $\Z_{q^2}$ coset contains at most $q$ elements of $S$, every $\Z_{pq^2}$ coset contains at most $q$ elements of $S$.

It follows that every $\Z_{p^2q^2}$ coset contains at most $pq$ elements of $S$, and hence
$$
|S|\leq pqr.
$$
Since $pqr\mid |S|$, we have $|S|=pqr$.

Therefore, every $\Z_{pq^2}$ coset contains exactly one $\Z_{q^2}$ coset having exactly $q$ elements of $S$. Since $q^2\nmid |S|$, we cannot have both $\Phi_q\mid m_{\La}$ and $\Phi_{q^2}\mid m_{\La}$. Thus, all nonempty intersections with $\Z_{q^2}$ cosets are either complete $\Z_q$ cosets or complete systems of representatives of the $\Z_q$ cosets. In either case, $S$ is a tile.

\item
Suppose that
$$
\Phi_{pr}\nmid m_{\La},\qquad
\Phi_{pqr}\nmid m_{\La},\qquad
\Phi_{pq^2r}\nmid m_{\La}.
$$
After translating, assume that $0\in S$. The conditions imply that, in every $\Z_{pq^2r}$ coset, the intersection with $S$ is contained either in a $\Z_{pq^2}$ coset or in a $\Z_{q^2r}$ coset.

In the first case, the intersection contains at most $pq$ elements. In the second case, it contains at most $qr$ elements. Let $l\leq p$ be the number of $\Z_{pq^2r}$ cosets in which $S$ is contained in a $\Z_{pq^2}$ coset. Then
$$
|S|\leq lpq+(p-l)qr.
$$
Since $|S|=kpqr$ for some $k$ with $\gcd(k,pqr)=1$, and since $p<r$, we obtain
$$
kpqr\leq lpq+(p-l)qr\leq pqr.
$$
Therefore, $k=1$ and $l=0$. Equality also implies that every $\Z_{pq^2r}$ coset contains exactly $qr$ elements of $S$, all contained in one $\Z_{q^2r}$ coset.

Fix one of these selected $\Z_{q^2r}$ cosets. It contains $r$ different $\Z_{q^2}$ cosets. Since $q^2\nmid |S|$, every $\Z_{q^2}$ coset contains at most $q$ elements of $S$. As the selected $\Z_{q^2r}$ coset contains exactly $qr$ elements of $S$, each of its $\Z_{q^2}$ cosets contains exactly $q$ elements of $S$.

We cannot have both $\Phi_q\mid m_{\La}$ and $\Phi_{q^2}\mid m_{\La}$, since this would imply $q^2\mid |\La|=|S|$. It follows that one of the following two possibilities holds in every $\Z_{q^2}$ coset that intersects $S$: the intersection is a complete $\Z_q$ coset, or it is a complete system of representatives of the $\Z_q$ cosets. The same possibility must hold in all the nonempty $\Z_{q^2}$ cosets. Otherwise, one intersection would give $\Phi_q\mid m_{\La}$ and another one would give $\Phi_{q^2}\mid m_{\La}$.

In either case, there is a common tiling complement $T_q$ in $\Z_{q^2}$. If all the intersections are systems of representatives, we take $T_q=\Z_q$. If all the intersections are $\Z_q$ cosets, we take $T_q$ to be any complete system of representatives of the $\Z_q$ cosets in $\Z_{q^2}$.

It follows that $S+T_q$ covers exactly once one complete $\Z_{q^2r}$ coset inside each $\Z_{pq^2r}$ coset. The selected $\Z_{q^2r}$ cosets form a complete system of representatives of the $\Z_p$ cosets in the quotient $\Z_{p^2q^2r}/\Z_{q^2r}\cong\Z_{p^2}$. Therefore,
$$
(S+T_q)+\Z_p=\Z_{p^2q^2r}.
$$
All the representations in this sum are unique. Hence,
$$
S+(T_q+\Z_p)=\Z_{p^2q^2r},
$$
and $S$ is a tile.

\item
Suppose that
$$
\Phi_{p^2}\nmid m_{\La},\qquad
\Phi_{p^2q}\nmid m_{\La},\qquad
\Phi_{p^2q^2}\nmid m_{\La}.
$$
The proof is similar to Case \ref{item:subcase1}. The set $S$ intersects every $\Z_{p^2q^2}$ coset in at most one $\Z_{pq^2}$ coset. Since each $\Z_{q^2}$ coset contains at most $q$ elements of $S$, every $\Z_{pq^2}$ coset, and hence every $\Z_{p^2q^2}$ coset, contains at most $pq$ elements of $S$.

Since $pqr\mid |S|$, it follows that $|S|=pqr$. Every $\Z_{p^2q^2}$ coset therefore intersects $S$ in exactly one $\Z_{pq^2}$ coset, and $S$ has exactly $pq$ elements in that coset. More precisely, it has exactly $q$ elements in every $\Z_{q^2}$ coset contained in that $\Z_{pq^2}$ coset.

Since $q^2\nmid |S|$, the nonempty intersections with the $\Z_{q^2}$ cosets are either all $\Z_q$ cosets or all systems of representatives of the $\Z_q$ cosets. In both cases, these intersections have a common tiling complement $T$ in $\Z_{q^2}$.

Therefore, $S'=S+T$ covers exactly once every $\Z_{pq^2}$ coset that intersects $S$. The set $S'$ contains exactly one $\Z_{pq^2}$ coset from each $\Z_{p^2q^2}$ coset. These selected $\Z_{pq^2}$ cosets have a common tiling complement $T'$ inside the corresponding $\Z_{p^2q^2}$ cosets. Since $S'$ intersects every $\Z_{p^2q^2}$ coset, we obtain
$$
\Z_{p^2q^2r}=S'+T'=(S+T)+T'=S+(T+T').
$$
Thus, $S$ tiles $\Z_{p^2q^2r}$.

\item
Suppose that
$$
\Phi_{p^2r}\nmid m_{\La},\qquad
\Phi_{p^2qr}\nmid m_{\La},\qquad
\Phi_{p^2q^2r}\nmid m_{\La}.
$$
If $\La$ is contained in a proper subgroup, then Lemma \ref{lem:red}.3 gives the result. Otherwise, these assumptions contradict Lemma \ref{lem:2primedonotdivide}, together with the corollary following it.
\end{enumerate}
\end{itemize}

\section{Appendix}

The main purpose of this Appendix is to prove Proposition \ref{prop:big}.
  This proposition provides a uniform treatment of the large-cardinality cases discussed in Remark \ref{remark:1}.  We first prove a statement about subsets of cyclic groups of prime-power order.

\begin{lem}\label{lem:padic}
Let $A\subseteq\Z_{p^n}$, and define
$$
\mathcal{D}(A)
=
\left\{
j\in\{1,\ldots,n\}:
\text{there are }a,a'\in A\text{ such that }a-a'
\text{ has order }p^j
\right\}.
$$
Then $|A|\leq p^{|\mathcal{D}(A)|}$. 
More generally, let $J\subseteq\{1,\ldots,n\}$ and suppose that
$$
\mathcal{D}(A)\subseteq J
\qquad\text{and}\qquad
|A|=p^{|J|}.
$$
Then $A$ tiles $\Z_{p^n}$. Moreover, a tiling complement can be chosen to depend only on $J$.
\end{lem}

\begin{proof}
Identify the elements of $\Z_{p^n}$ with the integers in
$\{0,1,\ldots,p^n-1\}$, and write
$$
a=a_0+a_1p+\cdots+a_{n-1}p^{n-1},
\qquad
a_i\in\{0,\ldots,p-1\}.
$$
If $a\neq a'$ and $i$ is the smallest index such that $a_i\neq a_i'$, then
$$
p^i\mid a-a',
\qquad
p^{i+1}\nmid a-a',
$$
and hence $a-a'$ has order $p^{n-i}$ in $\Z_{p^n}$.

Consider the rooted $p$-adic tree whose vertices at level $i$ are the residue classes modulo $p^i$ that contain an element of $A$. A vertex at level $i$ has at most $p$ children at level $i+1$. If some vertex at level $i$ has at least two children, then there are two elements of $A$ that agree in their first $i$ base-$p$ digits and differ in the next digit. Their difference has order $p^{n-i}$. Thus, branching at level $i$ can occur only if
$$
n-i\in\mathcal{D}(A).
$$
At every other level, each occupied vertex has only one occupied child. Therefore, the number of leaves is at most multiplied by $p$ at each of the $|\mathcal{D}(A)|$ possible branching levels. Since the leaves are the elements of $A$, we obtain
$$
|A|\leq p^{|\mathcal{D}(A)|}.
$$

Now suppose that $\mathcal{D}(A)\subseteq J$ and $|A|=p^{|J|}$. The first part gives
$$
p^{|J|}
=
|A|
\leq
p^{|\mathcal{D}(A)|}
\leq
p^{|J|}.
$$
Consequently,
$$
\mathcal{D}(A)=J.
$$

Define $I_J=\{n-j:j\in J\}\subseteq\{0,\ldots,n-1\}$ 
and
$$
T_J
=
\left\{
\sum_{i=0}^{n-1}t_ip^i:
t_i=0\text{ if }i\in I_J,\quad
t_i\in\{0,\ldots,p-1\}\text{ if }i\notin I_J
\right\}.
$$
Then $|T_J|=p^{n-|J|}$. 
Moreover, every nonzero difference of two elements of $T_J$ has order $p^j$ for some $j\notin J$.

Suppose that
 $a+t=a'+t'$ 
 for some $a,a'\in A$ and $t,t'\in T_J$. Then
$a-a'=t'-t$. 
If this difference is nonzero, its order belongs to $\mathcal{D}(A)=J$. On the other hand, since it is also a nonzero difference of two elements of $T_J$, its order is $p^j$ for some $j\notin J$. This is impossible. Therefore,
$a=a'$ and 
$t=t'$. 
Thus, all the representations in $A+T_J$ are unique. Finally,
$$
|A||T_J|
=
  p^{|J|}p^{n-|J|}=
p^n.
$$
It follows that
$$
A+T_J=\Z_{p^n}
$$
with unique representations. Hence, $A$ tiles $\Z_{p^n}$, and the complement $T_J$ depends only on $J$.
\end{proof}

\begin{prp}\label{prop:big}
Let $G$ be a finite cyclic group, and let $(S,\La)$ be a spectral pair in $G$. Suppose that
$$
\frac{|G|}{\gcd(|G|,|S|)}=p^k
$$
for some prime $p$ and some $k\in\N$. Then $S$ is a tile.
\end{prp}

\begin{proof}
Write $|G|=p^nd$, 
where $n\geq k$ and $p\nmid d$. The assumption gives
$$
\gcd(|G|,|S|)=p^{n-k}d.
$$
In particular,
$p^{n-k}d\mid |S|$ 
and $v_p(|S|)=n-k$. 
Using the Chinese remainder theorem, identify
$$
G\cong\Z_{p^n}\times\Z_d.
$$
For every $y\in\Z_d$, define the fiber
$$
S_y
=
\{x\in\Z_{p^n}:(x,y)\in S\}.
$$
There are $d$ such fibers. Since
$$
|S|\geq p^{n-k}d,
$$
the pigeonhole principle shows that there is a $y_0\in\Z_d$ such that
$$
|S_{y_0}|\geq p^{n-k}.
$$

Let
$$
J
=
\left\{
j\in\{1,\ldots,n\}:
\Phi_{p^j}\mid m_{\La}
\right\}.
$$
If $x,x'\in S_y$ and $x-x'$ has order $p^j$ in $\Z_{p^n}$, then the corresponding elements $(x,y),(x',y)\in S$ also have difference of order $p^j$ in $G$. Since $(S,\La)$ is a spectral pair, this implies
$\Phi_{p^j}\mid m_{\La}$.
Therefore, $
\mathcal{D}(S_y)\subseteq J$
for every $y\in\Z_d$.

Applying Lemma \ref{lem:padic} to the fiber $S_{y_0}$ gives
$$
p^{n-k}
\leq
|S_{y_0}|
\leq
p^{|\mathcal{D}(S_{y_0})|}
\leq
p^{|J|}.
$$
Hence, $|J|\geq n-k$.  
On the other hand, the cyclotomic polynomials $\Phi_{p^j}$, $j\in J$, are relatively prime in $\Z[x]$. Therefore,
$$
\prod_{j\in J}\Phi_{p^j}\mid m_{\La}.
$$
Evaluating at $x=1$ and using $\Phi_{p^j}(1)=p$, we obtain
$$
p^{|J|}
\mid
m_{\La}(1)
=|\La|=
|S|.
$$
Since $v_p(|S|)=n-k$, it follows that
$|J|\leq n-k$, 
consequently,
$|J|=n-k$. 

For every $y\in\Z_d$, Lemma \ref{lem:padic} now gives
$$
|S_y|
\leq
p^{|\mathcal{D}(S_y)|}
\leq
p^{|J|}
=
p^{n-k}.
$$
Therefore,
$$
|S|
=
\sum_{y\in\Z_d}|S_y|
\leq
dp^{n-k}.
$$
Since $p^{n-k}d\mid |S|$ and $S\neq\emptyset$, we conclude that
$$
|S|=p^{n-k}d.
$$
Equality in the preceding estimate implies that
$$
|S_y|=p^{n-k}=p^{|J|}
$$
for every $y\in\Z_d$.

Since $\mathcal{D}(S_y)\subseteq J$, the second part of Lemma \ref{lem:padic} shows that every fiber $S_y$ tiles $\Z_{p^n}$ with the same tiling complement $T_J$. In particular,
$$
S_y+T_J=\Z_{p^n}
$$
with unique representations for every $y\in\Z_d$.

Consider  $T_J$ as the subset $T_J\times\{0\}
\subseteq
\Z_{p^n}\times\Z_d$.  
For each fixed $y\in\Z_d$, we have
$$
(S_y\times\{y\})+(T_J\times\{0\})
=
\Z_{p^n}\times\{y\}
$$
with unique representations. Taking the union over all $y\in\Z_d$, we obtain
$$
S+(T_J\times\{0\})
=
\Z_{p^n}\times\Z_d
\cong
G
$$
with unique representations. Hence, $S$ is a tile of $G$.
\end{proof}

\section*{Acknowledgments}

We thank Izabella \L aba for her helpful comments and suggestions, which improved the paper. G.K. and G.S. also thank the Graduate Center of the City University of New York for its academic support during the completion of this work.

\printbibliography 

\bigskip

\noindent
\textbf{Thomas Fallon}\\
Brooklyn College, City University of New York\\
\href{mailto:fallon.t@gmail.com}{fallon.t@gmail.com}

\medskip

\noindent
\textbf{Gergely Kiss}\\
Corvinus University of Budapest\\
HUN-REN Alfr\'ed R\'enyi Institute of Mathematics, Analysis Department\\
\href{mailto:kigergo57@gmail.com}{kigergo57@gmail.com}

\medskip

\noindent
\textbf{Azita Mayeli}\\
The Graduate Center, City University of New York\\
\href{mailto:amayeli@gc.cuny.edu}{amayeli@gc.cuny.edu}

\medskip

\noindent
\textbf{G\'abor Somlai}\\
The University of Melbourne and
E\"otv\"os Lor\'and University\\
The Graduate Center, City University of New York\\
\href{mailto:gabor.somlai@ttk.elte.hu}{gabor.somlai@unimelb.edu.au}

\end{document}